\documentclass[a4paper, 12pt]{article}

\usepackage[sort&compress]{natbib}
\bibpunct{(}{)}{;}{a}{}{,} 

\usepackage{amsthm, amsmath, amssymb, mathrsfs, multirow, url, subfigure}
\usepackage{graphicx} 
\usepackage{ifthen} 
\usepackage{amsfonts}
\usepackage[usenames]{color}
\usepackage{fullpage}

\theoremstyle{plain}

\newtheorem*{vp}{Validity Principle}
\newtheorem*{ep}{Efficiency Principle}

\theoremstyle{definition}

\theoremstyle{remark} 

\newtheorem*{astep}{A--step}
\newtheorem*{pstep}{P--step}
\newtheorem*{cstep}{C--step}

\newcommand{\bel}{\mathsf{bel}}
\newcommand{\pl}{\mathsf{pl}}

\newcommand{\unif}{{\sf Unif}}
\newcommand{\nm}{{\sf N}}

\newcommand{\chisq}{{\sf ChiSq}}

\newcommand{\XX}{\mathbb{X}}

\newcommand{\UU}{\mathbb{U}}

\newcommand{\Xbar}{\bar X}
\newcommand{\xbar}{\overline{x}}

\renewcommand{\S}{\mathcal{S}}
\renewcommand{\SS}{\mathbb{S}}

\renewcommand{\phi}{\varphi}

\title{Frameworks for prior-free posterior probabilistic inference}
\author{
Chuanhai Liu \\
Department of Statistics \\
Purdue University \\
\url{chuanhai@purdue.edu} \\
\mbox{} \\
Ryan Martin \\
Department of Mathematics, Statistics, and Computer Science \\
University of Illinois at Chicago \\
\url{rgmartin@uic.edu} \\
}
\date{\today}

\begin{document}

\maketitle 

\begin{abstract}
The development of statistical methods for valid and efficient probabilistic inference without prior distributions has a long history.  Fisher's fiducial inference is perhaps the most famous of these attempts.  We argue that, despite its seemingly prior-free formulation, fiducial and its various extensions are not prior-free and, therefore, do not meet the requirements for prior-free probabilistic inference.  In contrast, the inferential model (IM) framework is genuinely prior-free and is shown to be a promising new method for generating both valid and efficient probabilistic inference.  With a brief introduction to the two fundamental principles, namely, the validity and efficiency principles, the three-step construction of the basic IM framework is discussed in the context of the validity principle.  Efficient IM methods, based on conditioning and marginalization are illustrated with two benchmark examples, namely, the bivariate normal with unknown correlation coefficient and the Behrens--Fisher problem.

\smallskip

\end{abstract}

\section{Introduction}
\label{S:intro}

Statistical inference is the process of converting experience, in the form of observed data, to knowledge about the underlying population in question, and is an essential part of the scientific method of discovery.  Our starting point in this paper will be a sampling model for observable data $X$, depending on a parameter $\theta$ in $\Theta$.  Mathematically, a sampling model for $X$ is a $\theta$-dependent probability distribution $P_{X|\theta}$ defined on the sample space $\XX$ of $X$.  The sampling model's dependence on $\theta$ implies that the observed data $X=x$ carries relevant information about the unknown parameter.  Our goal is to convert this information into a probabilistic measure of uncertainty.  That is, for an assertion or hypothesis $A \subseteq \Theta$ concerning the unknown parameter, we want to assign a measure of the ``plausibility'' that the assertion $A$ is true.  This plausibility measure should depend on the data $x$, should have a meaningful probabilistic interpretation, and should not require specification of an \emph{a priori} distribution for $\theta$.  

R.~A.~Fisher, starting with his work on inverse probability \citep{fisher1930}, made ambitious efforts to develop prior-free probabilistic inference based on a fiducial argument. \citet[][p.~54]{fisher1973} writes: 
\begin{quote}
By contrast [to the Bayesian argument], the fiducial argument uses the observations only to change the logical status of the parameter from one in which nothing is known of it, and no probability statement about it can be made, to the status of a random variable having a well-defined distribution.
\end{quote}
Although some of the fiducial ideas were reinterpreted by \citet{neyman1941} and used to create confidence intervals, a central concept in the frequentist paradigm, fiducial inference has been perceived as Fisher's ``one great failure'' \citep{zabell1992}.  Fisher acknowledged his only limited success in developing a framework for prior-free probabilistic inference based on the fiducial argument, he insisted that there was something valuable in it.  He wrote: 
\begin{quote}
I don't understand yet what fiducial probability does. We shall have to live with it a long time before we know what it's doing for us. But it should not be ignored just because we don't yet have a clear interpretation \citep[][p.~926]{savage1964}.
\end{quote}
Fisher's confidence in the value of fiducial inference has inspired continued efforts, including structural inference \citep{fraser1968}, Dempster--Shafer theory \citep{dempster2008, shafer1976}, generalized p-values and confidence intervals \citep{chiang2001, weerahandi1993}, generalized fiducial \citep{hannig2009, hannig2012, hannig.lee.2009}, confidence distributions \citep{xie.singh.strawderman.2011, xie.singh.2012}, and Bayesian inference with default, reference, and/or data-dependent priors \citep{berger2006, bergerbernardosun2009, fraser.reid.marras.yi.2010, fraser2011, mghosh2011}.   Here we will argue, however, that fiducial inference and its variants mentioned above, are actually not prior-free.

\citet{imbasics} have recently introduced an alternative paradigm under the name \emph{inferential models} (IMs); see, also, the References and Further Reading below.  With the focus on logical reasoning with uncertainty, these authors seek the best possible approach to genuinely prior-free probabilistic inference.  The focus here is mainly on its two fundamental principles, namely, the \emph{Validity} and \emph{Efficiency} Principles, and formal approaches built upon these two principles.  This includes basic IMs, conditional IMs for combining information, and marginal IMs for efficient inference on interest parameters.

\section{Difficulties in reasoning toward prior-free inference}

\subsection{Belief functions and models for total ignorance}

To represent realistic knowledge about assertions of interest, the use of lower and upper probabilities or belief and plausibility functions is necessary.  Given data $X$, they are functions from the space of all subsets of $\Theta$ to $[0, 1]$, denoted by 
\[ \bel_X(A) \quad \text{and} \quad \pl_X(A) = 1-\bel_X(A^c), \quad A \subseteq \Theta,\]
where $A^c = \Theta \setminus A$ stands for the negation of the assertion $A$.  The belief function $\bel_X(A)$ represents the amount of evidence in data $X$ supporting the claim that the assertion $A$ is true.  Thus, the value of $\bel_X(A^c)$ represents the amount of evidence
in data $X$ supporting the claim that the assertion $A^c$ is true or, equivalently, that the assertion $A$ is false.  The plausibility function $\pl_X(A)$ represents the amount of evidence in data $X$ that does not support the claim that $A$ is false.  It follows that $\bel_X(A) \leq \pl_X(A)$ for all assertions $A$.  Also, unlike conventional probabilities, which are additive, belief functions are sub-additive:
\begin{equation}
\label{eq:sub-additive}
\bel_X(A) + \bel_X(A^c) \leq 1, \quad \text{for all assertions $A$}. 
\end{equation}
For more discussion on belief functions, see \citet{shafer1976}.

In terms of belief functions, a model of total ignorance can be formally written as
\begin{equation}
\label{eq:total-ignorance}
\bel(A) = 0 \quad \text{and} \quad  \pl(A)=1 \quad \text{for all $A \subseteq \Theta$}.
\end{equation}
That is, no evidence is available to either support or refute any assertion $A$.  This definition is consistent with what is given in \citet[][p.~36]{fisher1973}: {\it ``The necessary {\em ignorance} is specified by our inability to discriminate any of the different sub-aggregates having different frequency ratios, such as must always exists.''}  It follows immediately that probability is not a satisfactory model for total ignorance and, therefore, no Bayesian priors can represent total ignorance; see \citet{walley1991}.  In other words, all Bayesian priors are informative, as the relative importance of any two points in $\Theta$ is fully specified.  Nevertheless, Bayesian methods based on conventional probabilities (with proper or improper priors) can be useful for constructing methods, e.g., confidence intervals, with good frequentist properties \citep{fraser2011}.

\subsection{Fiducial inference is not prior-free}

There is a close connection between fiducial and default/``non-informative'' prior Bayes, so the limitations of the latter, discussed above, must also be limitations of the former.  To see this, suppose that there is a joint distribution for $(X,\theta)$ that is consistent with both the sampling model (``$X | \theta$'') and the fiducial distribution (``$\theta | X$'').  Then $\theta$, or some transformation thereof, must be a location parameter, and the fiducial distribution corresponds to the Bayesian posterior obtained from a flat prior on the location parameter.  See \citet{lindley1958} and \citet{taraldsen.lindqvist.2013} for details.  

To make these points clear, we take a closer look at the fiducial argument or, more precisely, the fiducial operation in a simple example.  Consider inference about the mean of the unit normal model $\nm(\theta, 1)$ from a single observation $X$.  A simple association for this sampling model can be written as
\[ X = \theta + Z, \quad Z \sim \nm(0,1). \]
In the fiducial literature, $Z$ is called the pivotal quantity since its distribution is free of parameters.  When $X$ is observed, the fiducial argument continues to regard $Z$, now written as $Z=X-\theta$, as a $\nm(0,1)$-distributed random variable.  This distribution on $Z$, and its functional connection to $\theta$, admits a fiducial distribution for $\theta$:
\[ \theta \mid X \sim \nm(X, 1). \]
The ``continue to regard'' \citep{dempster1963} reasoning is what drives the operation that changes the knowledge status $\theta$ from total ignorance to that which can be represented by conventional probability.  This reasoning does not appear to be consistent with the goal of prior-free inference.  Indeed, the standard view that $\theta$ is fixed but unknown agrees with the model of total ignorance.  However, once $X$ is fixed at its observed value $x$, the conditional distribution of $Z$ degenerates at the point $x-\theta$.  The fiducial operation that replaces the degenerate conditional distribution $Z \sim \delta_{x-\theta}$ with the non-degenerate conditional distribution $Z \sim \nm(0,1)$ must be using some information beyond that contained in $x$ and $P_{X|\theta}$.  Therefore, fiducial inference cannot be prior-free.

\section{The IM framework}

\subsection{Background on the development}

The precise formulation of the IM framework will be given below, but it may be of interest to know how we reached this particular formulation since some might consider the journey to be at least as important as the destination.   The starting point was a precise statement of the goal of statistical inference: 
\begin{quote}
to give, for any assertion or hypothesis $A$ about the parameter of interest, meaningful summaries of the evidence in the observed data supporting the claims that ``$A$ is true'' and ``$A^c$ is true.''  
\end{quote}
Given these constraints, the goal is to make ``the best possible inference.''  This general idea motivates the two principles given in the next section.  

Towards this goal, a first question is if probability is the appropriate kind of measure.  In cases where a genuine prior distribution is available, the Bayesian approach based on posterior probabilities is appropriate.  However, when no meaningful prior information is available, as is often the case in scientific applications, difficulties arise.  First, as discussed above, there is no prior distribution that encodes ignorance.  Second, the use of default or non-informative priors that depend on the sampling model differs enough from the conventional Bayesian approach that, in our opinion, it cannot really be considered Bayes.  So, from a foundational point of view, we conclude that the Bayesian approach is not appropriate for prior-free probabilistic inference.  

The next question is if fiducial inference, or one of its variants, is appropriate.  As discussed above, fiducial has some difficulties.  In fact, fiducial is basically Bayes with a possibly data-dependent prior \citep{hannig2012}, so the same issues discussed above would apply to fiducial.  Dempster--Shafer theory was another candidate method, with the appeal of not basing inference on probabilities but on belief functions.  The difficulty with the Dempster--Shafer approach, in our opinion, is that the corresponding belief function values cannot be interpreted on a common scale.  That is, in one application, 0.9 might be a large belief function value, but in another it might be small.  Having a common scale on which the belief function values can be interpreted---this is the ``meaningfulness'' part of the above definition---is, in our opinion, essential the the success of a method as a tool for scientific discovery.  

So, after considering a variety of existing methods, Dempster--Shafer's use of belief functions was desirable, but a tool to properly calibrate the belief function values was needed.  The idea to expand the range between the belief and plausibility function values with a random set calibrated to the distribution of unobservable auxiliary variable was designed precisely to meet this need.  This approach achieves what is set out in the above definition; see, also, the Validity Principle below and \eqref{eq:bel.valid}.  The goal then to do the ``best possible,'' subject to the validity constraint, leads to the Efficiency Principle and optimality considerations \citep{imbasics}.    

Incidentally, though the IM belief function output is not a probability measure, we consider IM inference to be ``probabilistic'' in a certain sense.  The point is that it can be explained through the prediction of an unobservable quantity with a random set, and belief is just a probability with respect to the distribution of that random set.  This is discussed more in the following section.

\subsection{Two fundamental principles and the basic IM framework}

Philosophically, the IM framework for statistical inference is built on the following validity principle \citep{sts.discuss.2014}.

\begin{vp}
Probabilistic inference requires associating an unobservable but predictable quantity with the observable data and unknown parameter. Probabilities to be used for inference are obtained by valid prediction of the predictable quantity.
\end{vp}

To make the notion of ``predictable quantity'' precise, we consider an alternative description of the sampling model.  Specify an auxiliary variable $U$ in $\UU$ with distribution $P_U$ and a function $a$ such that the sample $X$, defined as 
\begin{equation}
\label{eq:sampling-model01}
X = a(\theta, U) \quad \text{where} \quad U \sim P_U, 
\end{equation}
has distribution $P_{X|\theta}$.  The pair $(a,P_U)$ is completely known.  

Unlike fiducial and its extensions, for valid prediction, the IM approach carries out its fundamental operations in the well-defined probability space of the auxiliary variable.  The key to this approach is the use of \emph{predictive random sets}.  According to \citet{imbasics}, a valid predictive random set $\S$ can be defined by specifying
\begin{itemize}
\item[(i)] a collection $\SS$ of nested subsets of $\UU$, including $\varnothing$ and $\UU$, to serve as the support of $\S$, and 
\vspace{-2mm}
\item[(ii)] the so-called natural measure $P_\S$ for $\S$ that satisfies
\begin{equation}\label{eq:PRS-measure}
P_\S(\S \subseteq K) = \sup_{S \in \SS: S \subseteq K} P_U(S), \quad K \subseteq \UU. 
\end{equation}
\end{itemize}
The centered random interval $\S[-|U|,|U|]$, where $U \sim \nm(0,1)$,
provides a simple example of a valid predictive random set for predicting a realization from $N(0,1)$.

The IM framework makes probabilistic inference by propagating prediction in the auxiliary variable space to the parameter space.  The three-step construction is as follows.  

\begin{astep}
Associate the observable data $X$ with the unknown quantity $\theta$ using an auxiliary variable $U \sim P_U$, such as that in \eqref{eq:sampling-model01}, to obtain the sets 
\begin{equation}
\label{eq:association-fromU2T}
\Theta_x(u) = \left\{\theta: x=a(u, \theta)\right\}, \quad x \in \XX, \quad u \in \UU.
\end{equation}
\end{astep}

\begin{pstep}
Predict $U$ using a valid predictive random set $\S$.
\end{pstep}

\begin{cstep}
Combine the association \eqref{eq:association-fromU2T} and prediction $\S$ to get
\begin{equation}
\label{eq:focal}
\Theta_x(\S) = \bigcup_{u \in \S} \Theta_x(u). 
\end{equation}
Then compute the belief and plausibility functions via the distribution of $\Theta_x(\S)$:
\begin{align} 
\bel_x(A;\S) & = P_\S\{\Theta_x(\S) \subseteq A\} \label{eq:basic-bel} \\
\pl_x(A;\S) & = P_\S\{\Theta_x(\S) \cap A \neq \varnothing\} \label{eq:basic-pl}
\end{align}
If $\Theta_x(\S) = \varnothing$, then $\S$ needs to be stretched, i.e., replace $\S$ with the smallest $S=S_\S$ in $\SS$ such that $S \supseteq \S$ and $\bigcup_{u \in S} \Theta_x(u)$ is non-empty \citep{leafliu2012}. 
\end{cstep}

The IM results are probabilistic and, since the predictive random set is valid, the IM results have desirable frequency properties \citep{imbasics}.  In particular, if $\S$ is valid, then for all $A \subseteq \Theta$ and all $\alpha \in (0,1)$, 
\begin{equation}
\label{eq:bel.valid}
\begin{split}
\sup_{\theta \not\in A} P_{X|\theta}\{\bel_X(A; \S) \geq 1-\alpha\} & \leq \alpha, \quad \text{and} \\ \sup_{\theta \in A} P_{X|\theta}\{\pl_X(A;\S) \leq \alpha\} & \leq \alpha.
\end{split}
\end{equation}
In other words, for example, $\pl_X(A;\S)$ is stochastically no smaller than $\unif(0,1)$, as a function of $X$, when the assertion $A$ is true.  This provides a meaningful and objective scale on which to interpret the plausibility (and belief) function values.  As a consequence, one can obtain procedures, such as tests and confidence regions, having exact control on frequentist error rates \citep{imbasics}.  

Further developments concern efficiency, motivated by the following principle.

\begin{ep}
Subject to the validity constraint, probabilistic inference should be made as efficient as possible.
\end{ep}

\citet{imcond, immarg} studied two classes of efficiency problems, namely, conditional IMs for combining information and marginal IMs for inference on interest parameters.  These are introduced below with two benchmark examples, which help demonstrate the differences between the IM and fiducial-type frameworks.

Finally, in some applications, prediction of future observations is the goal, not inference on $\theta$.  Prediction problems can be posed as ones involving marginalization or model averaging.  The IM framework is capable of handling such problems too, but this will not be discussed any further here; see \citet{impred} for details.

\subsection{Conditional IMs}

Take the bivariate normal model with zero means, unit variances, and unknown correlation coefficient $\theta$.  Consider inference about $\theta$ from a sample of size $n$, $\{(Y_{1,i}, Y_{2,i}):\; i=1,...,n\}.$  The data-generating based association, for example, has a $2n$-dimensional auxiliary variable.  \citet{imcond} argue that the fully observed functions of the auxiliary variable do not need to be predicted and, by conditioning on the observed components, the prediction of the unobserved components can be sharpened.  This argument implies that the association can be built from the model's sufficient statistics
\[ X_1 = \frac{1}{2}\sum_{i=1}^n (Y_{1,i}+Y_{2,i})^2 \quad \text{and} \quad X_2 = \frac{1}{2}\sum_{i=1}^n (Y_{1,i}-Y_{2,i})^2. \]
With independent chi-square auxiliary variables, $U_1$ and $U_2$, with $n$ degrees of freedom, the association becomes
\[ X_1 = (1+\theta)U_1 \quad \text{and} \quad X_2 = (1-\theta)U_2. \]
It is tempting to follow a fiducial-type argument and condition on $X_1/U_1 + X_2/U_2 = 2$, an attractive parameter-free identity.  However, such conditioning amounts to conditioning on all the data $(X_1, X_2)$, which makes the predictive distribution degenerate and, consequently, the corresponding inference is not valid for all assertions.

Note that the auxiliary variable is two-dimensional while the parameter of interest is one-dimensional.  To obtain an association with a lower-dimensional auxiliary variable, \citet{imcond} propose to first rewrite the above association in terms of a new pair of auxiliary variables $(V_1,V_2)$ such that only $V_1$, say, carries information directly related to $\theta$.  By doing so, a conditional association obtains by taking $V_1$ as the auxiliary variable, with distribution conditioned on the observed value of $V_2$.  This leads to a reduction of the auxiliary dimension from two to one.  For the bivariate normal correlation problem, their approach based on a partial differential equations gives the conditioning function or component of the form
\[ V_2 \equiv (1+\theta)\log U_1 + (1-\theta) \log U_2. \]
Note that the function above depends on $\theta$; \citet{imcond} add the adjective ``local'' in such cases.  Inference can then be made by predicting a corresponding complement transformation of $(U_1, U_2)$, for example, 
\[ V_1 \equiv \ln U_1 - \ln U_2. \]
Let $F_\theta$ denote the distribution function of $V_1$ conditional on 
\[ V_2 = (1+\theta)\log \frac{X_1}{1+\theta} + (1-\theta) \log \frac{X_2}{1-\theta}. \]
The so-called local conditional IM can be represented using the basic IM with the auxiliary variable $W$ and the association
\[ \log\frac{X_1}{X_2} = \log\frac{1+\theta}{1-\theta} + F^{-1}_\theta(W) \quad W \sim \unif(0,1).
\]
For more technical details, see \citet{imcond}.  The numerical results presented there show that this local conditional IM produces exact confidence intervals that are more efficient than the best of available frequentist methods.

\subsection{Marginal IMs}

Suppose that $\theta$ is a multidimensional parameter, but that only some lower-dimensional feature of $\theta$ is of interest.  If $\theta$ itself were the parameter of interest, then we could find an association to connect it to data and a set of auxiliary variables.  However, in the case where only a feature of $\theta$ is of interest, efficiency demands that we further reduce the dimension of the auxiliary variable.  The marginal IM approach of \citet{immarg} presents a strategy for efficient, and often optimal inference on interest parameters.  We choose the benchmark Behrens--Fisher problem to illustrate MIMs and to contrast with the fiducial approach.  

Suppose we have two independent samples of size $n_1$ and $n_2$ from $\nm(\mu_1, \sigma_1^2)$ and $\nm(\mu_2, \sigma_2^2)$, respectively, where all four normal parameters are unknown.  The Behrens--Fisher problem concerns marginal inference on the difference $\mu_2-\mu_1$.  The CIM theory suggests that we start with an association based on the sampling distribution of the sufficient statistics, i.e., sample means and variances $\bar{X}_k$ and $S^2_k$, $k=1,2$.
Let $\delta=\mu_2-\mu_1$ be the parameter of interest.  Take the following association, in terms of $\delta$ and the nuisance parameters $\mu_1$, $\sigma_1^2$ and $\sigma_2^2$,
\begin{equation}
\label{eq:BF-association-delta}
\bar{Y}\equiv \bar{X}_2-\bar{X}_1 = \delta + \bigl(\sigma_1^2/n_1+\sigma_2^2/n_2\bigr)^{1/2}
Z, 
\end{equation}
\begin{equation}
\label{eq:BF-association-nuisance}
\bar{X}_1 = \mu_1 + \sigma_1 n_1^{-1/2} U, \quad S_k^2 = \sigma_k^2 V_k^2\qquad k=1,2,
\end{equation}
where $V_1$ and $V_2$ are independent with $V_k^2 \sim \chisq(n_k-1)/(n_{k}-1)$ for $k=1,2$, $(Z, U)$ is bivariate normal with zero means, unit variances, and correlation coefficient depending on $\sigma_k$'s and $n_k$'s, and $(V_1,V_2)$ and $(Z,U)$ are independent.
Fisher's fiducial solution is obtained by ``continuing to regard'' the auxiliary variable as having its sampling distribution, conditional on the observed sufficient statistics.  From the IM perspective, this fiducial argument is questionable.  It is also easy to check that this fiducial distribution is the same as the Bayesian posterior distribution based on assigning flat priors to the ``location parameters'' $\mu_1$, $\mu_2$, $\log\sigma_1^2$, and $\log \sigma_2^2$.  

Valid inference can be made with any admissible predictive random set for the four-dimensional auxiliary variable $(Z, U, V_1^2, V_2^2)$.  The shape of the predictive random set controls the precision on features of the auxiliary variable.  Due to the unknown $\mu_1$ in the first expression of \eqref{eq:BF-association-nuisance}, the accuracy in predicting $U$ is not useful because it provides no information on $\sigma_1$ and, therefore, results in less efficient inference on $\delta$ in \eqref{eq:BF-association-delta}.  This suggests that we stretch the predictive random set along the $U$-coordinate as much as possible, ending up with an effectively three-dimensional or marginalized predictive random set for $(Z, V_1^2, V_2^2)$.  This logical reasoning for efficient inference leads to the initial dimension-reduced association:
\begin{align*}
\bar{Y} & = \delta + \left(\sigma_1^2/n_1+\sigma_2^2/n_2\right)^{1/2} Z, \\ 
S_k^2 & = \sigma_k^2 V_k^2 \qquad k=1,2, 
\end{align*}
with the same $(Z, V_1^2, V_2^2)$.  That same argument for the usefulness of accurate prediction of $\sigma_1^2$ for inferring $\delta$ leads to the further dimension-reduced association
\begin{equation}
\label{eq:BF-2D}
\frac{\bar{Y}-\delta}{\sqrt{S_1^2/n_1}}
= (1+\lambda^2 )^{1/2}\frac{Z}{V_1} \quad \text{and} \quad \frac{S_2^2/n_2}{S_1^2/n_1} = \lambda^2\frac{V_2^2}{V_1^2}
\end{equation}
for inference about $\delta$, with $\lambda^2 = (n_1\sigma_2^2)/(n_2\sigma_1^2)$ as the new nuisance parameter.  Thus, without loss of efficiency, we can make valid inference on $\delta$ by predicting a two-dimensional auxiliary variable $(Z/V_1, V_2^2/V_1^2)$.  However, it does not appear that the same argument can be used further to obtain an association for $\delta$ involving only a one-dimensional auxiliary variable.

To associate $\delta$ to a one-dimensional auxiliary variable, \citet{immarg} allows the auxiliary variable to depend on nuisance parameters.  Such a nuisance parameter-dependent auxiliary variable is constructed in such a way that the effect of the nuisance parameter on its distribution is minimized in some sense.  A heuristic approach to constructing such a $\lambda^2$-dependent auxiliary variable from \eqref{eq:BF-2D} is as follows.  For values of $V_2^2/V_1^2$ such that $V_2^2/V_1^2 \approx 1$, the value of $1+\lambda^2 V_2^2/V_1^2$ is approximately $1+\lambda^2$.
This suggests we divide the first identity of \eqref{eq:BF-2D}
by
\[ \bigl(1+(S_2^2/n_2)/(S_1^2/n_1)\bigr)^{1/2} = \bigl( 1+ \lambda^2V_2^2/V_1^2 \bigr)^{1/2}, \]
a variant of the second expression in \eqref{eq:BF-2D}.  This gives
\begin{equation}
\label{eq:BF-tomega}
\frac{\bar{Y}-\delta}{\sqrt{\frac{S_1^2}{n_1}+\frac{S_2^2}{n_2}}} = U_\omega \equiv \frac{Z}{\sqrt{\omega V_1^2+ (1-\omega)V_2^2}},
\end{equation}
where $\omega = 1/(1+\lambda^2)$.  Thus inference about $\delta$ can be made by a predictive random set for the one-dimensional quantity $U_\omega$ alone.  Let $F_\omega$ denote the distribution function of $U_\omega$.  We rewrite the above association as
\begin{equation}
\label{eq:BF-1D}
\frac{\bar{Y}-\delta}{\sqrt{S_1^2/n_1+S_2^2/n_2}} = F_\omega^{-1}(U), \quad U \sim \unif(0,1).
\end{equation}

Introduce the predictive random set $\S=\{u:\; |u-0.5| \leq |U'-0.5|\}$, where $U' \sim \unif(0,1)$, for $U$ in \eqref{eq:BF-1D}.  The corresponding predictive random set for $U_\omega$ in \eqref{eq:BF-tomega} is
\[ F^{-1}_\omega(\S) \equiv \{F_\omega^{-1}(u):\; u \in \S\}. \]
To eliminate the effect of $\omega$, take the enlarged predictive random set
\[ \S^+ = \bigcup_\omega F^{-1}_\omega(\S). \]
This predictive random set is valid for predicting $U_\omega$ in \eqref{eq:BF-tomega}, uniformly in $\omega$ and, therefore, leads to valid inference on assertions about $\delta$.  It can be shown that $\S^+$ is equivalent to the centered predictive random set for predicting a realization from the t-distribution with $\min\{n_1-1, n_2-1\}$ degrees of freedom \citep{immarg}.

\subsection{An example}

Suppose that $X_1,\ldots,X_n$ are iid $\nm(\mu,\sigma^2)$.  Here, both $\mu$ and $\sigma^2$ are unknown, but, for the moment, the parameter of interest is $\psi=\mu/\sigma$, the standardized mean, or signal-to-noise ratio.  We may first reduce to the joint minimal sufficient statistics $(\Xbar, S^2)$, the sample mean and sample variance, respectively.  Based on this pair, the (conditional) association can be written as 
\begin{align*}
\Xbar & = \mu + \sigma n^{-1/2} U_1, \quad & U_1 & \sim \nm(0,1), \\
S & = \sigma U_2, \quad & (n-1)U_2^2 & \sim \chisq(n-1),
\end{align*}
where $U_1$ and $U_2$ are independent.  The parameter of interest is a scalar but the auxiliary variable $(U_1,U_2)$ is two-dimensional, so we would like to further reduce the dimension of the latter.  Then the conditional association can be rewritten as 
\[ \frac{n^{1/2}\Xbar}{S} = \frac{n^{1/2}\psi + U_1}{U_2}, \quad \text{and} \quad S = \sigma U_2. \]
The first expression has no dependence on $\sigma$.  Also, for any pair $(S,U_2)$, there exists a $\sigma$ that satisfies the second expression.  Therefore, in the language of \citet{immarg}, this is a ``regular'' association so the nuisance parameter $\sigma$ can be eliminated by ignoring the second expression above.  That is, the marginal association is 
\[ \frac{n^{1/2}\Xbar}{S} = \frac{n^{1/2}\psi + U_1}{U_2}. \]
The right-hand side of the above expression has a known distribution, namely, a non-central Student-t distribution with $n-1$ degrees of freedom and non-centrality parameter $n^{1/2}\psi$.  Therefore, a simple change of auxiliary variable gives the modified marginal association:
\[ \frac{n^{1/2}\Xbar}{S} = F_{n,\psi}^{-1}(W), \quad W \sim \unif(0,1), \]
where $F_{n,\psi}$ is the above non-central Student-t distribution function.  To this point, our calculations agree with those of \citet{dempster1963} using fiducial arguments.  

To complete the marginal IM for $\psi$, the A-step gives the sets 
\[ \Psi_x(w) = \{\psi: F_{n,\psi}(t_x) = w \}, \quad w \in (0,1), \]
where $t_x = n^{1/2}\xbar/s$.  For the P-step, we propose to use a simple ``default'' \citep{imbasics} predictive random set:
\[ \S = \{w \in (0,1): |w-0.5| \leq |W-0.5|\}, \quad W \sim \unif(0,1). \]
This predictive random set satisfies the required properties for the corresponding IM to be valid \citep{imbasics}.  Finally, the C-step gives 
\[ \Psi_x(\S) \equiv \bigcup_{u \in \S} \Psi_x(u) = \{\psi: | F_{n,\psi}(t_x) - 0.5 | \leq | W - 0.5| \}. \]
Since the predictive random set is valid, one can produce valid probabilistic inference for any assertion about $\psi$ of interest.  In particular, for singleton assertions, i.e., $\{\psi\}$, the belief function is zero but the marginal plausibility function is 
\[ \pl_x(\psi; \S) = 1 - |2 F_{n,\psi}(t_x) - 1 |. \]
We may construct an interval estimate for $\psi$ using this marginal plausibility function.  Specifically, for some fixed level $\alpha \in (0,1)$, the $100(1-\alpha)$\% plausibility interval for $\psi$ is $\{\psi: \pl_x(\psi; \S) > \alpha\}$ which, in this case, simplifies to 
\[ \{\psi: \pl_x(\psi; \S) > \alpha\} = \{\psi: \alpha / 2 < F_{n,\psi}(t_x) < 1-\alpha/2 \}. \]
This plausibility interval clearly has frequentist coverage probability $1-\alpha$, which is a consequence of our general IM validity results.  

\begin{figure}[t]
\label{fig:cv}
\begin{center}
\subfigure[$\mu=1$]{\scalebox{0.45}{\includegraphics{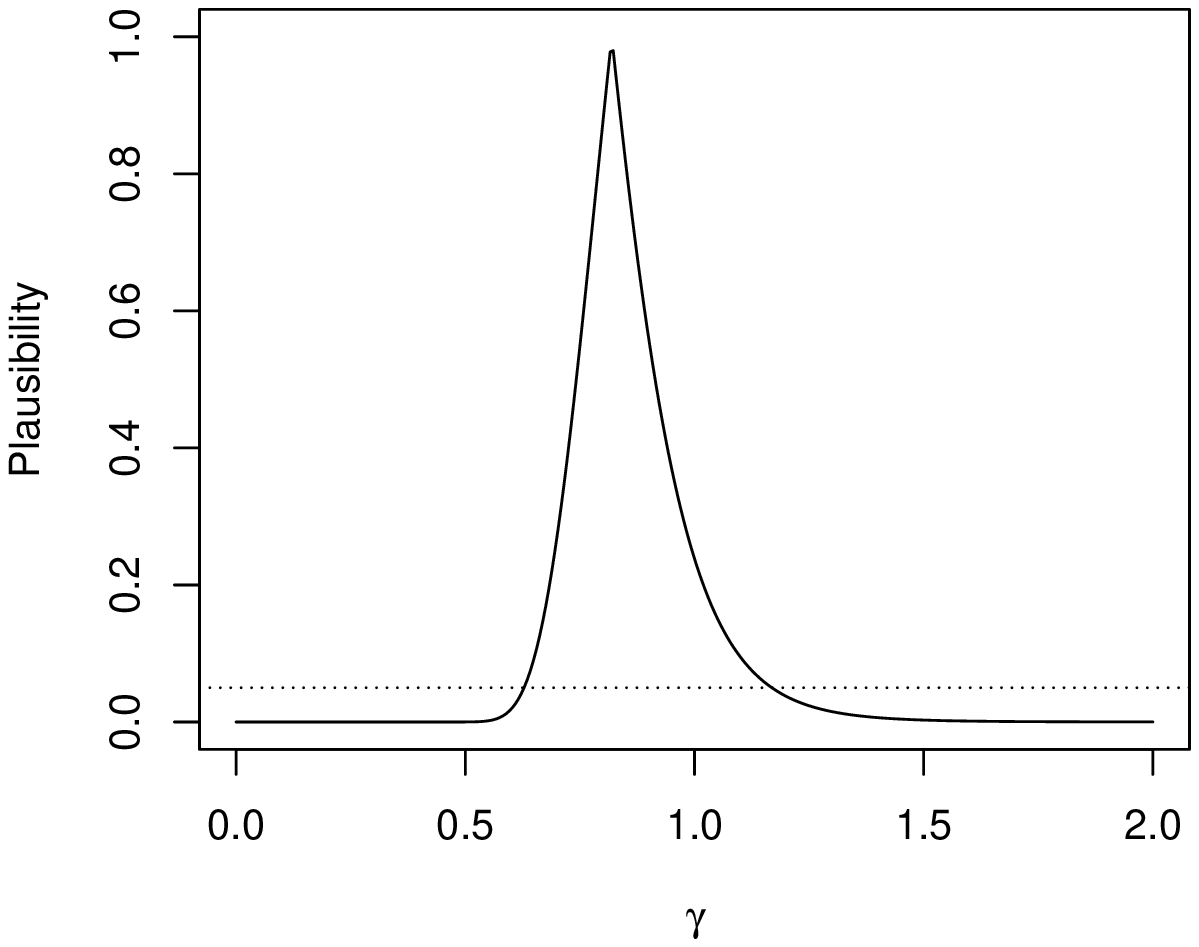}}} 
\subfigure[$\mu=0$]{\scalebox{0.45}{\includegraphics{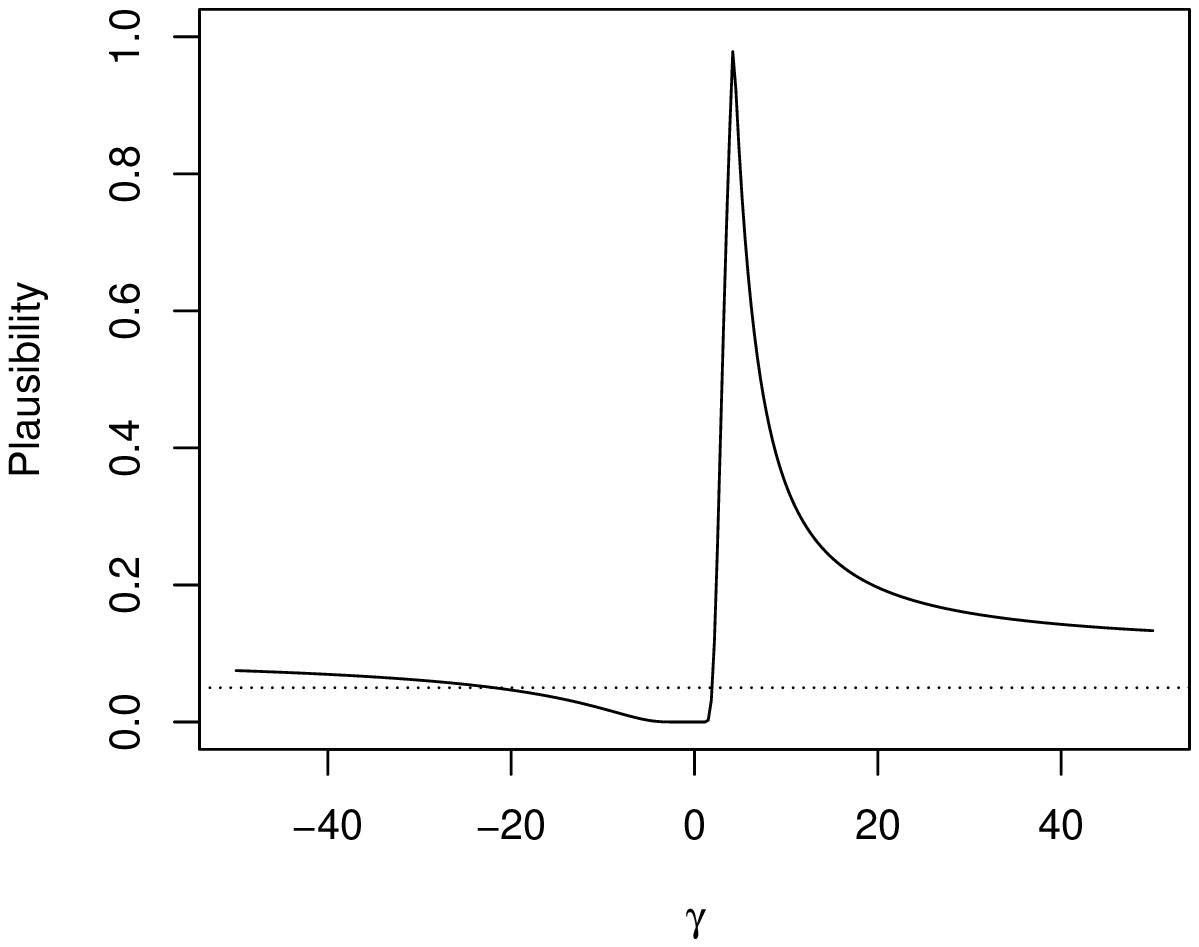}}}
\end{center}
\caption{Plots of the plausibility function for $\gamma=\sigma/\mu$, the coefficient of variation, based on samples of size $n=50$ from $\nm(\mu,1)$ for two values of $\mu$.}
\end{figure}

As a variation on this example, next consider $\gamma = 1/\psi = \sigma/\mu$, the coefficient of variation.  \citet{berger.liseo.wolpert.1999} demonstrate that this is a challenging problem for Bayesian, likelihood, and frequentist frameworks.  Using the approach just described for the signal-to-noise ratio problem, it is straightforward to construct a marginal IM for $\gamma$.  In particular, for singleton assertions, the plausibility function is 
\[ \pl_x(\gamma; \S) = 1 - |2 F_{n,1/\gamma}(t_x) - 1 |. \]
Plots of this plausibility function for two samples of size $n=50$ from $\nm(\mu,1)$, for two values of $\mu$, are displayed in Figure~\ref{fig:cv}.  Panel~(a) shows the case of $\mu=1$, and it is clear that the data are fairly informative about $\gamma$, and the true $\gamma=1$ is contained in the 95\% plausibility interval.  Panel~(b) shows the case of $\mu=0$.  In this case, data are not particularly informative, i.e., the plausibility function does not seem to vanish as $\gamma \to \pm\infty$, so the plausibility regions are unbounded.  Since we know the IM plausibility regions are exact, it follows from the theorem of \citet{gleser.hwang.1987} that they must be unbounded with positive probability.  However, we argue that this is not a problem of the IM approach.  Given the usual variation of the sample mean around $\mu$, if $\mu$ is (close to) zero, then it is not possible to rule out $\gamma$ values with very large magnitude.  So, in such cases, no reasonable approach should be able to rule out these extreme values of $\gamma$ based on data alone.  Therefore, the marginal IM approach here arguably gives the ``best possible'' inference on $\gamma$.

\section{Conclusion}
\label{S:discuss}

Methods for prior-free probabilistic inference are fundamentally important for converting experience to knowledge in scientific explorations.  This is especially true in the ``big data'' world we now live in.  In spite of its failure as a general method leading to prior-free probabilistic inference, Fisher's fiducial arguments have inspired many statisticians to think outside the box and to develop promising new approaches.  For example, J.~Hannig and his coauthors have developed a generalized fiducial framework that yields inferential methods which share the asymptotic efficiency of classical likelihood-based methods but often perform significantly better in small-sample studies.  

The IM framework is motivated by Fisher's fiducial argument, but the two differ both philosophically and technically.  This new approach is promising in that it is truly prior-free and produces posterior probabilistic assessments of uncertainty with desirable frequency properties.  This is really a new school of thought, with the ambitious goal to carry out the best possible sampling model-based scientific inference, so the efforts so far have focused on developing the fundamental ideas, building blocks of the framework.  Additional work can be found below in the References and Further Reading sections.  As discussed briefly here, the current work on conditional and marginal IMs shows that the development of efficient IMs will be very interesting in years to come.  While there is still so much to do concerning theory, computation, and application of IMs, the authors believe that it has a very bright future given the fundamental role that statistics has to play in the development of science.

\section*{Acknowledgments}

The authors thank the Editor, Associate Editor, and reviewers for helpful suggestions on an earlier draft of this paper.  This research is partially supported by the U.S.~National Science Foundation, grants DMS--1007678, DMS--1208833, and DMS--1208841.



\end{document}